\documentclass[11pt]{article}
\usepackage{amssymb}
\usepackage{amsmath}
\usepackage[all]{xy}
\usepackage{times}
\usepackage{graphicx}

\title{Borel chromatic number of closed graphs\indent}
\author{Dominique LECOMTE and Miroslav ZELEN\'Y$^1$}
\date{\today}

\def\ufootnote#1{\let\savedthfn\thefootnote\let\thefootnote\relax
\footnote{#1}\let\thefootnote\savedthfn\addtocounter{footnote}{-1}}

\newcommand{\boraxi}{{\bf\Sigma}^{0}_{\xi}}

\newcommand{\bortwo}{{\bf\Delta}^{0}_{2}}

\newcommand{\bormone}{{\bf\Pi}^{0}_{1}}

\newcommand{\bormtwo}{{\bf\Pi}^{0}_{2}}

\newcommand{\bormxi}{{\bf\Pi}^{0}_{\xi}}

\newcommand{\borxi}{{\bf\Delta}^{0}_{\xi}}

\newtheorem{thm} {Theorem} [section]

\newtheorem{lem} [thm] {Lemma}

\newtheorem{defi} [thm] {Definition}

\def\Baire{\omega^\omega}

\begin{document}

\maketitle

\centerline{$\bullet$ Universit\' e Paris 6, Institut de Math\'ematiques de Jussieu, Projet Analyse Fonctionnelle}

\centerline{Couloir 16-26, 4\`eme \'etage, Case 247, 4, place Jussieu, 75 252 Paris Cedex 05, France}

\centerline{dominique.lecomte@upmc.fr}\bigskip

\centerline{$\bullet$ Universit\'e de Picardie, I.U.T. de l'Oise, site de Creil,}

\centerline{13, all\'ee de la fa\"\i encerie, 60 107 Creil, France}\bigskip

\centerline{$\bullet^1$ Charles University, Faculty of Mathematics and Physics, Department of Mathematical Analysis}

\centerline{Sokolovsk\'a 83, 186 75 Prague, Czech Republic}

\centerline{zeleny@karlin.mff.cuni.cz}\bigskip\bigskip\bigskip\bigskip\bigskip\bigskip

\ufootnote{{\it 2010 Mathematics Subject Classification.}~Primary: 03E15, Secondary: 54H05}

\ufootnote{{\it Keywords and phrases.}~Borel chromatic number, Borel class, coloring}

\ufootnote{{\it Acknowledgements.}~The main result was obtained during the first author's 
stay at Charles University in Prague in May 2014. The first author thanks Charles University in Prague for the hospitality.}

\ufootnote{The research was supported by the grant GA{\v C}R P201/12/0436 for the second author.}

\noindent {\bf Abstract.} We construct, for each countable ordinal $\xi$, a closed graph with Borel chromatic number two and Baire class $\xi$ chromatic number $\aleph_0$. 

\vfill\eject

\section{$\!\!\!\!\!\!$ Introduction}\indent

 The study of the Borel chromatic number of analytic graphs on Polish spaces was initiated in 
[K-S-T]. In particular, the authors prove in this paper that the Borel chromatic number of the graph generated by a partial Borel function has to be in $\{ 1,2,3,\aleph_0\}$. They also provide a minimum graph ${\cal G}_0$ of uncountable Borel chromatic number. This last result had a lot of developments. For example, B. Miller gave in [Mi] some other versions of it, which helped him to generalize a number of known dichotomy theorems in descriptive set theory. The first author generalized in [L2] the ${\cal G}_0$-dichotomy to any dimension making sense in classical descriptive set theory, and also used versions of ${\cal G}_0$ to study the non-potentially closed subsets of a product of two Polish spaces (see [L1]).\bigskip

 A study of the $\borxi$ chromatic number of analytic graphs on Polish spaces was initiated in [L-Z1] and was motivated by the ${\cal G}_0$-dichotomy. More precisely, let $B$ be a Borel binary relation, on a Polish space $X$, having a Borel countable coloring (i.e., a Borel map $c\! :\! X\!\rightarrow\!\omega$ such that $c(x)\!\not=\! c(y)$ if $(x,y)\!\in\! B$). Is there a relation between the Borel class of $B$ and that of the coloring?  In other words, is there a map 
$k\! :\!\omega_1\!\setminus\!\{ 0\}\!\rightarrow\!\omega_1\!\setminus\!\{ 0\}$ such that any $\bormxi$ binary relation having a Borel countable coloring has in fact a ${\bf\Delta}^0_{k(\xi )}$-measurable countable coloring, for each $\xi\!\in\!\omega_1\!\setminus\!\{ 0\}$?\bigskip

 In [L-Z2], the authors give a negative answer: for each countable ordinal $\xi\!\geq\! 1$, there is a partial injection with disjoint domain and range 
$i\! :\!\omega^\omega\!\rightarrow\!\omega^\omega$, whose graph\smallskip

\noindent - is $D_2(\bormone )$ (i.e., the difference of two closed sets),\smallskip

\noindent - has Borel chomatic number two,\smallskip

\noindent - has no $\borxi$-measurable countable coloring.\bigskip

 On the other hand, they note that an open binary relation having a finite coloring $c$ has also a 
$\bortwo$-measurable finite coloring (consider the differences of the 
$\overline{c^{-1}(\{ n\} )}$'s, for $n$ in the range of the coloring). Note that an irreflexive closed binary relation on a zero-dimensional space has a continuous countable coloring (this coloring is $\bortwo$-measurable in non zero-dimensional spaces). So they wonder whether we can build, for each countable ordinal $\xi\!\geq\! 1$, a closed binary relation with a Borel finite coloring but no 
$\borxi$-measurable finite coloring. This is indeed the case:\bigskip

\noindent\bf Theorem\it\ Let $\xi\!\geq\! 1$ be a countable ordinal. Then there exists a partial injection with disjoint domain and range $f\! :\!\omega^\omega\!\rightarrow\!\omega^\omega$ whose graph is closed (and thus has Borel chromatic number two), and has no $\borxi$-measurable finite coloring (and thus has $\borxi$ chromatic number $\aleph_0$).\rm\bigskip

 The previous discussion shows that this result is optimal. Its proof uses, among other things, the method used in [L-Z2] improving Theorem 4 in [M]. This method relates topological complexity and Baire category.

\section{$\!\!\!\!\!\!$ M\'atrai sets}\indent

 Before proving our main result, we recall some material from [L-Z2].\bigskip

\noindent\bf Notation.\rm\ The symbol $\tau$ denotes the usual product topology on the Baire space $\omega^\omega$. 

\begin{defi} We say that a partial map $f\! :\!\omega^\omega\!\rightarrow\!\omega^\omega$ is 
\bf nice\rm\ if its graph $\mbox{Gr}(f)$ is a $(\tau\!\times\!\tau )$-closed subset of 
$\omega^\omega\!\times\!\omega^\omega$.\end{defi}

 The construction of $P_\xi$ and $\tau_\xi$, and the verification of the properties (1)-(3) from the next lemma (a corollary of Lemma 2.6 in [L-Z2]), can be found in [M], up to minor modifications.
 
\begin{lem} \label{Matrai} Let $1 \leq \xi < \omega_1$. Then there are 
$P_\xi\!\subseteq\!\omega^\omega$, and a topology $\tau_\xi$ on $\omega^\omega$ such that\smallskip

(1) $\tau_\xi$ is zero-dimensional perfect Polish and 
$\tau\!\subseteq\!\tau_\xi\!\subseteq\!\boraxi (\tau)$,\smallskip

(2) $P_\xi$ is a nonempty $\tau_\xi$-closed nowhere dense set,\smallskip

(3) if $S\!\in\!\boraxi (\Baire ,\tau )$ is $\tau_\xi$-nonmeager in $P_\xi$, then $S$ is 
$\tau_\xi$-nonmeager in $\omega^\omega$,\smallskip

(4) if $V,W$ are nonempty $\tau_\xi$-open subsets of $\omega^\omega$, then we can find a 
$\tau_\xi$-dense $G_\delta$ subset $H$ of $V\!\setminus\! P_\xi$, a $\tau_\xi$-dense $G_\delta$ subset $L$ of $W\!\setminus\! P_\xi$, and a nice $(\tau_\xi ,\tau_\xi )$-homeomorphism from $H$ onto $L$.\end{lem}

 The following lemma (a corollary of Lemma 2.7 in [L-Z2]) is a consequence of the previous one. It provides, among other things, a topology $T_\xi$ that we will use in the sequel.

\begin{lem} \label{Txi} Let $1 \leq \xi < \omega_1$. Then there is a disjoint countable family 
${\cal G}_\xi$ of subsets of $\Baire$ and a topology $T_\xi$ on $\Baire$ such that\smallskip

(a) $T_\xi$ is zero-dimensional perfect Polish and 
$\tau\!\subseteq\! T_\xi\!\subseteq\!\boraxi (\tau )$,\smallskip

(b) for any nonempty $T_\xi$-open sets $V,V'$, there are disjoint $G,G'\!\in\! {\cal G}_\xi$ with 
$G\!\subseteq\! V$, $G'\!\subseteq\! V'$, and there is a nice $(T_\xi ,T_\xi )$-homeomorphism from $G$ onto $G'$,\smallskip

\noindent and, for every $G\!\in\! {\cal G}_\xi$,\smallskip

(c) $G$ is nonempty, $T_\xi$-nowhere dense, and in $\bormtwo (T_\xi )$,\smallskip

(d) if $S\!\in\!\boraxi (\Baire ,\tau )$ is $T_\xi$-nonmeager in $G$, then $S$ is $T_\xi$-nonmeager in 
$\Baire$.\end{lem}

 The construction of ${\cal G}_\xi$ and $T_\xi$ ensures that $T_\xi$ is $(\tau_\xi )^\omega$, where 
$\tau_\xi$ is as in Lemma \ref{Matrai}. This topology is on $(\omega^\omega )^\omega$, identified with 
$\omega^\omega$. We will need the following consequence of the construction of ${\cal G}_\xi$ and $T_\xi$.

\begin{lem} \label{compact} Let $1 \leq \xi < \omega_1$, and $V$ be a nonempty $T_\xi$-open set. Then $\overline{V}^\tau$ is not $\tau$-compact.\end{lem}

\noindent\bf Proof.\rm\ The fact that $T_\xi$ is $(\tau_\xi )^\omega$ gives a finite sequence $U_0, ...,U_n$ of nonempty open subsets of $(\omega^\omega ,\tau_\xi )$ with 
$U_0\!\times\! ...\!\times\! U_n\!\times\! (\omega^\omega )^\omega\!\subseteq\! V$. Thus 
$\overline{V}^\tau$ contains the $\tau$-closed set 
$\overline{U_0}^\tau\!\times\! ...\!\times\!\overline{U_n}^\tau\!\times\! (\omega^\omega )^\omega$, and it is enough to see that this last set is not $\tau$-compact. This comes from the fact that the Baire space $(\omega^\omega ,\tau )$ is not compact.\hfill{$\square$}

\section{$\!\!\!\!\!\!$ Proof of the main result}\indent

 Before proving our main result, we give an example giving the flavour of the sequel. In [Za], the author gives a Hurewicz-like test to see when two disjoint subsets $A,B$ of a product 
$Y\!\times\! Z$ of Polish spaces can be separated by an open rectangle. We set 
$\mathbb{A}\! :=\!\{ (n^\infty ,n^\infty )\mid n\!\in\!\omega\}$, 
$$\mathbb{B}_0\! :=\!\big\{\big( 0^{m+1}(n\! +\! 1)^\infty ,(m\! +\! 1)^{n+1}0^\infty\big)\mid 
m,n\!\in\!\omega\big\}$$ 
and $\mathbb{B}_1\! :=\!\big\{\big( (m\! +\! 1)^{n+1}0^\infty ,0^{m+1}(n\! +\! 1)^\infty\big)\mid 
m,n\!\in\!\omega\big\}$. Then $A$ is not separable from $B$ by an open rectangle exactly when 
there are $\varepsilon\!\in\! 2$ and continuous maps $g\! :\!\omega^\omega\!\rightarrow\! Y$, 
$h\! :\!\omega^\omega\!\rightarrow\! Z$ such that $\mathbb{A}\!\subseteq\! (g\!\times\! h)^{-1}(A)$ and $\mathbb{B}_\varepsilon\!\subseteq\! (g\!\times\! h)^{-1}(B)$.

\vfill\eject

\noindent\bf Example.\rm\ Here we are looking for closed graphs with Borel chromatic number two and of arbitrarily high finite $\borxi$ chromatic number $n$. There is an example with $\xi\! =\! 1$ and $n\! =\! 3$ where $\mathbb{B}_0$ is involved. We set 
$C\! :=\!\big\{\big( (2m)^\infty ,(2m\! +\! 1)^\infty\big)\mid m\!\in\!\omega\big\}\cup\mathbb{B}_0$,  
$$D\! :=\!\{ (2m)^\infty\mid m\!\in\!\omega\}\cup\{ 0^{m+1}(n\! +\! 1)^\infty\mid m,n\!\in\!\omega\}
\mbox{,}$$  
$R\! :=\!\{ (2m\! +\! 1)^\infty\mid m\!\in\!\omega\}\cup\{ (m\! +\! 1)^{n+1}0^\infty\mid m,n\!\in\!\omega\}$,  
$$f\big( (2m)^\infty\big)\! :=\! (2m\! +\! 1)^\infty\mbox{ and }
f\big( 0^{m+1}(n\! +\! 1)^\infty\big)\! :=\! (m\! +\! 1)^{n+1}0^\infty .$$ 
This defines $f\! :\! D\!\rightarrow\! R$ whose graph is $C$. The first part of $C$ is discrete, and thus closed. Assume that $(\alpha_k,\beta_k)\! :=\!
\big( 0^{m_k+1}(n_k\! +\! 1)^\infty ,(m_k\! +\! 1)^{n_k+1}0^\infty\big)\!\in\!\mathbb{B}_0$ and converges to $(\alpha ,\beta )\!\in\!\omega^\omega\!\times\!\omega^\omega$ as $k$ goes to infinity. We may assume that $(m_k)$ is constant, and $(n_k)$ too, so that 
$(\alpha ,\beta )\!\in\!\mathbb{B}_0$, which is therefore closed. This shows that $C$ is closed. Note that $D,R$ are disjoint and Borel, so that $C$ has Borel chromatic number two. Let $\Delta$ be a clopen subset of $\omega^\omega$. Let us prove that $C\cap\Delta^2$ or $C\cap (\neg\Delta )^2$ is not empty. We argue by contradiction. Then $\Delta$ or $\neg\Delta$ has to contain $0^\infty$. Assume that it is $\Delta$, the other case being similar. Then $0^{m+1}(n\! +\! 1)^\infty\!\in\!\Delta$ if $m$ is big enough. Thus $(m\! +\! 1)^{n+1}0^\infty\!\notin\!\Delta$ if $m$ is big enough. Therefore 
$(m\! +\! 1)^\infty\!\notin\!\Delta$ if $m$ is big enough. Thus 
$\big( (2m)^\infty ,(2m\! +\! 1)^\infty\big)\!\in\! C\cap (\neg\Delta )^2$ if $m$ is big enough, which is absurd.\bigskip

 We now turn to the general case. Our main lemma is as follows. We equip $\omega^m$ with the discrete topology $\tau_d$, for each $m\! >\! 0$.\bigskip

\noindent\bf Lemma\it\ Let $\xi\!\geq\! 1$ be a countable ordinal, $n\!\geq\! 1$ be a natural number, and $X\! :=\!\omega\!\times\!\omega^\omega$. Then we can find a partial injection 
$f\! :\! X\!\rightarrow\! X$ and a disjoint countable family $\cal F$ of subsets of $X$ such that\smallskip

(a) $f$ has disjoint domain and range,\smallskip

(b) $\mbox{Gr}(f)$ is $\big( (\tau_d\!\times\!\tau )\!\times\! (\tau_d\!\times\!\tau )\big)$-closed,\smallskip

(c) there is no sequence $(\Delta_i)_{i<n}$ of $\borxi$ subsets of $(X,\tau_d\!\times\!\tau )$ such that\smallskip

~~~~~~~(i) $\forall i\! <\! n~~\mbox{Gr}(f)\cap\Delta_i^2\! =\!\emptyset$,\smallskip
 
~~~~~~~(ii) $\bigcup_{i<n}~\Delta_i$ is $(\tau_d\!\times\! T_\xi )$-comeager in $X$,\smallskip

(d) $\cal F$ has the properties (b)-(d) in Lemma \ref{Txi}, where ${\cal G}_\xi$, $\omega^\omega$, 
$T_\xi$ and $\tau$ are respectively replaced with $\cal F$, $X$, $\tau_d\!\times\! T_\xi$ and 
$\tau_d\!\times\!\tau$,\smallskip
 
(e) $(\bigcup {\cal F})\cap\big(\mbox{Domain}(f)\cup\mbox{Range}(f)\big)\! =\!\emptyset$.\rm\bigskip

\noindent\bf Proof.\rm\  We argue by induction on $n$.\bigskip

\noindent\bf The basic case $n=1$\rm\bigskip

 Let ${\cal G}_\xi$ be the family given by Lemma \ref{Txi}. We split ${\cal G}_\xi$ into two disjoint subfamilies ${\cal G}^0_\xi$ and ${\cal G}^1_\xi$ having the property (b) in Lemma \ref{Txi}. This is possible since the elements of ${\cal G}_\xi$ are $T_\xi$-nowhere dense. Let 
$G_0,G_1\!\in\! {\cal G}^0_\xi$ be disjoint, and $\varphi$ be a nice $(T_\xi ,T_\xi )$-homeomorphism from $G_0$ onto $G_1$. We then set 
$f(0,\alpha )\! :=\!\big( 0,\varphi (\alpha )\big)$ if $\alpha\!\in\! G_0$, and 
${\cal F}\! :=\!\big\{\{ n\}\!\times\! G\mid n\!\in\!\omega\wedge G\!\in\! {\cal G}^1_\xi\big\}$. It remains to check that the property (c) is satisfied. We argue by contradiction, which gives 
$\Delta_0\!\in\!\borxi$. By property (d) in Lemma \ref{Txi}, 
$\Delta_0\cap (\{ 0\}\!\times\! G_\varepsilon )$ is $(\tau_d\!\times\! T_\xi )$-comeager in 
$\{ 0\}\!\times\! G_\varepsilon$ for each $\varepsilon\!\in\! 2$. As $f$ is a 
$(\tau_d\!\times\! T_\xi ,\tau_d\!\times\! T_\xi )$-homeomorphism, 
$\Delta_0\cap (\{ 0\}\!\times\! G_0)\cap f^{-1}\big(\Delta_0\cap (\{ 0\}\!\times\! G_1)\big)$ is 
$(\tau_d\!\times\! T_\xi )$-comeager in $\{ 0\}\!\times\! G_0$, which contradicts the fact that 
$\mbox{Gr}(f)\cap\Delta_0^2\! =\!\emptyset$.

\vfill\eject

\noindent\bf The induction step from $n$ to $n\! +\! 1$\rm\bigskip

 The induction assumption gives $f$ and $\cal F$. Here again, we split $\cal F$ into two disjoint subfamilies ${\cal F}^0$ and ${\cal F}^1$ having the property (b) in Lemma \ref{Txi}, where 
${\cal G}_\xi$, $\omega^\omega$, $T_\xi$ and $\tau$ are respectively replaced with $\cal F$, $X$, 
$\tau_d\!\times\! T_\xi$ and $\tau_d\!\times\!\tau$. Let $(V_p)$ be a basis for the topology 
$\tau_d\!\times\! T_\xi$ made of nonempty sets. Fix $p\!\in\!\omega$. By Lemma \ref{compact}, there is a countable family $(W^p_q)_{q\in\omega}$, with $(\tau_d\!\times\!\tau )$-closed union, and made of pairwise disjoint $(\tau_d\!\times\!\tau )$-clopen subsets of $X$ intersecting $V_p$.\bigskip

\noindent $\bullet$ Let $b\! :\!\omega\!\rightarrow\!\omega^2$ be a bijection. We construct, for 
$\vec v\! =\! (p,q)\!\in\!\omega^2$ and $\varepsilon\!\in\! 2$, and by induction on $b^{-1}(\vec v)$,\bigskip

- $G^{\vec v}_\varepsilon\!\in\! {\cal F}^0$,\smallskip

- a nice $(\tau_d\!\times\! T_\xi ,\tau_d\!\times\! T_\xi )$-homeomorphism 
$\varphi^{\vec v}\! :\! G^{\vec v}_0\!\rightarrow\! G^{\vec v}_1$.\bigskip

\noindent We want these objects to satisfy the following:\bigskip

- $G^{\vec v}_0\!\subseteq\! (V_p\cap W^p_q)\!\setminus\! 
(\bigcup_{m<b^{-1}(\vec v)}~\overline{G^{b(m)}_0\cup G^{b(m)}_1}^{\tau_d\times T_\xi})$,\bigskip

- $G^{\vec v}_1\!\subseteq\! V_q\!\setminus\! (G^{\vec v}_0\cup
\bigcup_{m<b^{-1}(\vec v)}~\overline{G^{b(m)}_0\cup G^{b(m)}_1}^{\tau_d\times T_\xi})$.\bigskip

\noindent $\bullet$ We now define the desired partial map $\tilde f\! :\!\omega\!\times\!\omega\!\times\!\omega^\omega\!\rightarrow\!\omega\!\times\!\omega\!\times\!\omega^\omega$, as well as  
$\tilde {\cal F}\!\subseteq\! 2^{\omega\times\omega\times\omega^\omega}$, as follows:
$$\tilde f(l,x)\! :=\!\left\{\!\!\!\!\!\!\!\!
\begin{array}{ll}
& \big( p\! +\! 1,\varphi^{p,q}(x)\big)\mbox{ if }l\! =\! 0~\wedge ~x\!\in\! G^{p,q}_0\mbox{,}\cr\cr
& \big( l,f(x)\big)\mbox{ if }l\! >\! 0~\wedge ~x\!\in\!\mbox{Domain}(f).
\end{array}
\right.$$
and $\tilde {\cal F}\! :=\!\big\{\{ l\}\!\times\! G\mid l\!\in\!\omega ~\wedge ~G\!\in\! {\cal F}^1\big\}$. Note that $\tilde f$ is well-defined and injective, by disjointness of the $(G^{\vec v}_0\cup G^{\vec v}_1)$'s. Identifying $X$ with $\omega\!\times\!\omega\!\times\!\omega^\omega$, we can consider $\tilde f$ as a partial map from $X$ into itself and $\tilde {\cal F}$ as a family of subsets of $X$ (this identification is based on the identification of $\omega$ with $\omega\!\times\!\omega$).\bigskip

\noindent (a), (d) and (e) are clearly satisfied.\bigskip

\noindent (b) Assume that $\big( (l_k,x_k),(m_k,y_k)\big)\!\in\!\mbox{Gr}(\tilde f)$ tends to 
$\big( (l,x),(m,y)\big)\!\in\! (\omega\!\times\! X)^2$ as $k$ goes to infinity. We may assume that 
$(l_k)$ and $(m_k)$ are constant.\bigskip

 If $l\! =\! 0$, then there is $p$ such that $p\! +\! 1\! =\! m$ and 
$(x_k,y_k)\!\in\! G^{p,q_k}_0\!\times\! G^{p,q_k}_1$. As $G^{p,q_k}_0\!\subseteq\! W^p_{q_k}$, we may also assume that $(q_k)$ is also constant and equals $q$. As $\varphi^{p,q}$ is nice, 
$\big( (l,x),(m,y)\big)\!\in\!\mbox{Gr}(\tilde f)$.\bigskip

 If $l\! >\! 0$, then $(x_k,y_k)\!\in\!\mbox{Gr}(f)$. As $\mbox{Gr}(f)$ is 
$\big( (\tau_d\!\times\!\tau )\!\times\! (\tau_d\!\times\!\tau )\big)$-closed, 
$\big( (l,x),(m,y)\big)\!\in\!\mbox{Gr}(\tilde f)$.\bigskip
 
\noindent (c) We argue by contradiction, which gives $(\Delta_i)_{i\leq n}$. We may assume, without loss of generality, that $(\{ 0\}\!\times\!\omega\!\times\!\omega^\omega )\cap\Delta_n$ 
is not meager in $(\{ 0\}\!\times\!\omega\!\times\!\omega^\omega ,\tau_d\!\times\! T_\xi )$. This gives 
$p\!\in\!\omega$ such that $(\{ 0\}\!\times\! V_p)\cap\Delta_n$ is $(\tau_d\!\times\! T_\xi )$-comeager in $V'_p\! :=\!\{ 0\}\!\times\! V_p$. As $V'_p\!\setminus\!\Delta_n\!\in\!\boraxi (\tau_d\!\times\!\tau )$, 
$(\{ 0\}\!\times\! G^{p,q}_0)\cap\Delta_n$ is $(\tau_d\!\times\! T_\xi )$-comeager in 
$\{ 0\}\!\times\! G^{p,q}_0$ for each $q\!\in\!\omega$.\bigskip

 As $\mbox{Gr}(\tilde f)\cap\Delta_n^2\! =\!\emptyset$ and the $\varphi^{\vec v}$'s are 
$(\tau_d\!\times\! T_\xi ,\tau_d\!\times\! T_\xi )$-homeomorphisms, 
$(\{ p+1\}\!\times\! G^{p,q}_1)\cap\Delta_n$ is $(\tau_d\!\times\! T_\xi )$-meager in 
$\{ p+1\}\!\times\! G^{p,q}_1$, for each $q$.\bigskip

 As $(\omega\!\times\!\omega\!\times\!\omega^\omega )\!\setminus\! (\bigcup_{i\leq n}~\Delta_i)$ is 
$(\tau_d\!\times\! T_\xi )$-meager in $\omega\!\times\!\omega\!\times\!\omega^\omega$ and 
$\borxi (\tau_d\!\times\!\tau )$, 
$$(\{ p+1\}\!\times\! G^{p,q}_1)\!\setminus\! (\bigcup_{i\leq n}~\Delta_i)$$ 
is $(\tau_d\!\times\! T_\xi )$-meager in $\{ p+1\}\!\times\! G^{p,q}_1$, for each $q$. Thus 
$(\{ p+1\}\!\times\! G^{p,q}_1)\cap (\bigcup_{i<n}~\Delta_i)$ is $(\tau_d\!\times\! T_\xi )$-comeager in 
$\{ p+1\}\!\times\! G^{p,q}_1$, for each $q$.\bigskip

\noindent\bf Claim \it\ The set 
$(\{ p+1\}\!\times\!\omega\!\times\!\omega^\omega )\cap (\bigcup_{i<n}~\Delta_i)$ is 
$(\tau_d\!\times\! T_\xi )$-comeager in $\{ p+1\}\!\times\!\omega\!\times\!\omega^\omega$.\rm\bigskip

 Indeed, we argue by contradiction. This gives $W\!\in\! (\tau_d\!\times\! T_\xi )\!\setminus\!\{\emptyset\}$ such that 
 $$(\{ p+1\}\!\times\! W)\cap (\bigcup_{i<n}~\Delta_i)$$ 
 is $(\tau_d\!\times\! T_\xi )$-meager in $W'\! :=\!\{ p+1\}\!\times\! W$. Let $q\!\in\!\omega$ be such that $V_q\!\subseteq\! W$. Then ${G^{p,q}_1\!\subseteq\! W}$ and 
$\{ p+1\}\!\times\! G^{p,q}_1\!\subseteq\! W'$. As 
$W'\cap (\bigcup_{i<n}~\Delta_i)\!\in\!\boraxi (\tau_d\!\times\!\tau )$ and 
$(\{ p+1\}\!\times\! G^{p,q}_1)\cap W'\cap (\bigcup_{i<n}~\Delta_i)$ is 
$(\tau_d\!\times\! T_\xi )$-comeager in $\{ p+1\}\!\times\! G^{p,q}_1$, 
$W'\cap (\bigcup_{i<n}~\Delta_i)$ is not $(\tau_d\!\times\! T_\xi )$-meager in $W'$, which is absurd.
\hfill{$\diamond$}\bigskip

 Now we set $\Delta'_i\! :=\! (\{ p+1\}\!\times\!\omega\!\times\!\omega^\omega )\cap\Delta_i$ if $i\! <\! n$. Note that $\Delta'_i\!\in\!\borxi (\{ p+1\}\!\times\!\omega\!\times\!\omega^\omega ,\tau_d\!\times\!\tau )$, 
$\mbox{Gr}(\tilde f)\cap (\Delta'_i )^2\! =\!\emptyset$, and $\bigcup_{i<n}~\Delta'_i$ is 
$(\tau_d\!\times\! T_\xi )$-comeager in $\{ p+1\}\!\times\!\omega\!\times\!\omega^\omega$, which contradicts the induction assumption.\hfill{$\square$}\bigskip

 In order to get our main result, it is enough to apply the main lemma to each $n\!\geq\! 1$. This gives 
$f_n\! :\!\omega\!\times\!\omega^\omega\!\rightarrow\!\omega\!\times\!\omega^\omega$. It remains to define $f\! :\!\bigcup_{n\geq 1}~(\{ n\}\!\times\!\omega\!\times\!\omega^\omega )\!\rightarrow\!
\bigcup_{n\geq 1}~(\{ n\}\!\times\!\omega\!\times\!\omega^\omega )$ by $f(n,x)\! :=\! f_n(x)$ (we identify 
$(\omega\!\setminus\!\{ 0\} )\!\times\!\omega\!\times\!\omega^\omega$ with $\omega^\omega$).

\section{$\!\!\!\!\!\!$ References}

\baselineskip=13.4pt

\noindent [K-S-T]\ \ A. S. Kechris, S. Solecki and S. Todor\v cevi\'c, Borel chromatic numbers,\ \it 
Adv. Math.\rm\ 141 (1999), 1-44

\noindent [L1]\ \ D. Lecomte, On minimal non potentially closed subsets of the plane,\ \it Topology Appl.\rm\ 154, 1 (2007), 241-262

\noindent [L2]\ \ D. Lecomte, A dichotomy characterizing analytic graphs of uncountable Borel chromatic number in any dimension,~\it Trans. Amer. Math. Soc.\rm~361 (2009), 4181-4193

\noindent [L-Z1]\ \ D. Lecomte and M. Zelen\'y, Baire-class $\xi$ colorings: the first three levels,\ \it Trans. Amer. Math. Soc.\rm\ 366, 5 (2014), 2345-2373

\noindent [L-Z2]\ \ D. Lecomte and M. Zelen\'y, Descriptive complexity of countable unions of Borel rectangles,\ \it Topology Appl.\rm\ 166 (2014), 66-84

\noindent [M]\ \ T. M\'atrai, On the closure of Baire classes under transfinite convergences,~\it Fund. Math.\ \rm 183, 2 (2004), 157-168

\noindent [Mi]\ \ B. Miller, The graph-theoretic approach to descriptive set theory,\ \it Bull. Symbolic Logic\rm\ 18, 4 (2012), 554-575 

\noindent [Za]\ \ R. Zamora, Separation of analytic sets by rectangles of low complexity,~\it manuscript (see arXiv)\ \rm
 
\end{document}